

\baselineskip=14pt
\parskip=10pt
\def\halmos{\hbox{\vrule height0.15cm width0.01cm\vbox{\hrule height
  0.01cm width0.2cm \vskip0.15cm \hrule height 0.01cm width0.2cm}\vrule
  height0.15cm width 0.01cm}}
\font\eightrm=cmr8 
\font\eighttt=cmtt8
\magnification=\magstephalf

\def\1{{\overline{1}}}
\def\2{{\overline{2}}}
\parindent=0pt
\overfullrule=0in

\def\frac#1#2{{#1 \over #2}}
\centerline
{
\bf
How to Extend K\'arolyi and  Nagy's BRILLIANT Proof of the Zeilberger-Bressoud 
}
\centerline
{
\bf
q-Dyson Theorem in order to Evaluate ANY Coefficient of the q-Dyson Product
}
\bigskip
\centerline{ {\it Shalosh B. EKHAD and
Doron 
ZEILBERGER}\footnote{$^1$}
{\eightrm  \raggedright
Department of Mathematics, Rutgers University (New Brunswick),
Hill Center-Busch Campus, 110 Frelinghuysen Rd., Piscataway,
NJ 08854-8019, USA.
{\eighttt zeilberg  at math dot rutgers dot edu} ,
\hfill \break
{\eighttt http://www.math.rutgers.edu/\~{}zeilberg/} .
Aug. 15, 2013.
Accompanied by Maple package \hfill \break {\eighttt qDYSON}
downloadable from Zeilberger's website.
Supported in part by the NSF.
}
}

\quad\quad\quad\quad\quad\quad\quad\quad\quad\quad\quad\quad
{\it Dedicated to Freeman Dyson (b. Dec. 15, 1923) on his $89 \frac{2}{3}$-th  birthday}

{\bf Very Important}: As in all our joint papers, this document, the human-readable {\it article},
is not the main point, but its Maple implementation, {\eighttt qDYSON}, written by DZ, available from

{\eighttt http://www.math.rutgers.edu/\~{}zeilberg/tokhniot/qDYSON} .

The `front' of this article

{\eighttt http://www.math.rutgers.edu/\~{}zeilberg/mamarim/mamarimhtml/qdyson.html} ,

contains sample input and output files, chuckfull of rigorously-derived deep identities,
computed by SBE.

In fact, since we believe in {\it free open access},
we supply the {\it source code}, that, in principle,
is also {\it humanly-readable}, but only if one knows Maple, and since most people do not know  it
fluently enough, DZ kindly prepared the present article. 

Let's first do a redux of Gyula K\'arolyi and Zolt\'an L\'or\'ant Nagy's {\it proof from the book}[KN],
of the Zeilberger-Bressoud[ZB] theorem (n\'ee Andrews's $q$-Dyson Conjecture [An]),
in a form that would be amenable for the extension promised in the title.

{\bf The K\'arolyi-Nagy  Brilliant Proof of Zeilberger-Bressoud}

{\bf Fact 1}: If a polynomial of degree $\leq d$ vanishes at $d+1$ {\bf different} places it must be identically zero.

{\bf Proof}: By induction on $d$. If $d=0$ it is a {\bf constant}, and since it happens to be zero
{\it somewhere}, it must be zero {\it everywhere}, i.e. it must be  the polynomial $0$. 

If a polynomial $P(x)$ of degree $\leq d$ vanishes at $d+1$ distinct places, let $a$ be one of them.
Then (thanks to Euclid), $P(x)=(x-a)Q(x)+c$ for some polynomial $Q(x)$ of degree $\leq d-1$ and constant $c$, that must be $0$, since $P(a)=0$.
Hence $P(x)=(x-a)Q(x)$, and $Q(x)$ is a polynomial of degree $\leq d-1$ that vanishes at $d$ different places.
By the induction hypothesis, it is the zero polynomial, and hence so is $P(x)$. \halmos

{\bf Fact 2}: (Lagrange Interpolation Formula) 
If $P(x)$ is a polynomial of degree $\leq  d$ in $x$, and $|A|=d+1$, then
$$
P(x)= 
\sum_{c \in A}^{} \left ( \, 
\prod_{c' \in A \backslash c}^{} \frac{x-c'}{c-c'} \, \right ) P(c) \quad .
$$
{\bf Proof}: The left-side minus the right side is a polynomial of degree $\leq d$ that vanishes at the $d+1$ different 
members of $A$, and hence must be identically zero by Fact 1. \halmos

{\bf Fact 3}: (Immediate consequence of the Lagrange Interpolation Formula) 
If $P(x)$ is a polynomial of degree $\leq  d$ in $x$, and $|A|=d+1$, then
$$
Coeff_{x^d} P(x)= 
\sum_{c \in A}^{} \frac{P(c)}
{\prod_{c' \in A \backslash c}^{} (c-c') }  \quad .
$$

{\bf Proof:} Extract the coefficient of $x^d$ on both sides of Fact 2. \halmos

{\bf Fact 4}: (`Quantitative' form ([L][KP]) of the Combinatorial Nullstellenatz ([A]), reproved in [KN], Lemma 2.1)
Let $F(x_1, \dots, x_n)$ be a polynomial of degree $\leq d_1 + \dots + d_n$. For arbitrary sets
$A_1, \dots, A_n$ with $|A_i|=d_i+1$, the coefficient of $\prod_{i=1}^{n} x_i^{d_i}$ in $F(x_1, \dots, x_n)$ is
$$
\sum_{c_1 \in A_1}  \sum_{c_2 \in A_2}  \dots \sum_{c_n \in A_n}  
\frac{F(c_1, \dots, c_n)}{\phi'_1(c_1)\phi'_2(c_2) \cdots \phi'_n(c_n)} \quad ,
$$
where $\phi_i(z):=\prod_{a \in A_i} (z-a)$, and $\phi'_i(c_i)$ is $\phi_i(z)/(z-c_i)$ evaluated at $c_i$.

{\bf New Proof}: Since any {\bf poly}nomial of  of degree $\leq d_1 + \dots + d_n$ is a linear
combination of {\bf mono}mials of degree $\leq d_1 + \dots + d_n$, it suffices, by {\bf linearity}, to prove this
for  monomials $F=x_1^{m_1} \cdots x_n^{m_n}$ with $m_1 + \dots + m_n \leq d_1 + \dots +d_n$.
If there is
an $i$ such that $m_i<d_i$ then the left side it $0$, and the right side is $0$ by Fact $3$ with
$x=x_i$, $d=d_i$, $P(x_i)=x_i^{m_i}$. If all $m_i \geq d_i$, then, of course, $m_i=d_i$ for all $1 \leq i \leq n$,
and the left side is $1$ and the right side is $1^n=1$, by applying Fact 3 for each $i$, with $P(x_i)=x_i^{d_i}$
and multiplying. \halmos

We are now ready for

{\bf Fact 5}: ({\bf The Zeilberger-Bressoud Theorem}([ZB])) Let, $q$ and $X$ be commuting indeterminates, and let $n$ be
a non-negative integer. Define first
$$
(X)_n :=\prod_{t=0}^{n-1} (1-q^t X) \quad .
$$
Let $a_1, \dots, a_n$ be non-negative integers, and let $x_1, \dots, x_n$ be commuting indeterminates.
The coefficient of $x_1^0 \dots x_n^0$ (i.e. the {\bf constant term}) of
$$
\prod_{1 \leq i < j \leq n} (x_i/x_j)_{a_i} (qx_j/x_i)_{a_j} \quad
\eqno(qDyson)
$$
equals the $q$-{\bf multinomial coefficient}
$$
\frac{(q)_{a_1+a_2+ \dots + a_n}}{(q)_{a_1} (q)_{a_2} \cdots (q)_{a_n}} \quad .
$$

{\bf Proof} ([KN] with purely-routine stuff removed).

If any of the $a_i$ equals $0$ then the theorem reduces to one with $<n$ variables and would follow
by induction on $n$, hence we can assume that all the $a_i$ are strictly positive.

Let $\sigma=\sum_{i=1}^{n} a_i$. We have to evaluate the coefficient of
$\prod_{i=1}^{n} x_i^{\sigma -a_i}$ of the polynomial (of degree $(n-1)\sigma$) 
$$
F(x_1, \dots, x_n):=\prod_{1 \leq i < j \leq n}  (x_i/x_j)_{a_i} (qx_j/x_i)_{a_j}  \, \cdot \, x_j^{a_i} x_i^{a_j} \quad .
$$
Let's Apply Fact 4 with $F$, $d_i=\sigma -a_i$ and, for $i=1, \dots, n$,
$$
A_i:=\{q^{\alpha_i} \, ; \, 0 \leq \alpha_i \leq \sigma -a_i \} \quad.
$$

{\bf SubFact 5.1}: If there exists a pair $i,j$ with  $1 \leq i<j \leq n$ such that
$$
-(a_i-1) \leq \alpha_i - \alpha_j \leq a_j  \quad ,
$$
then $F(q^{\alpha_1}, \dots , q^{\alpha_n})=0$. 

{\bf Proof}: Since $(q^{-e})_{f}=0$  if $0 \leq e <f$, 
 $(x_i/x_j)_{a_i} (qx_j/x_i)_{a_j}$ with $x_i=q^{\alpha_i}$ and $x_j=q^{\alpha_j}$ vanishes, and hence so does $F$.
\halmos

{\bf SubFact 5.2}: The set of lattice points $(\alpha_1, \dots, \alpha_n)$ with
$0 \leq \alpha_i \leq \sigma - a_i$ such that for {\it every} pair $1 \leq i <j \leq n$
it is not the case that $-(a_i-1) \leq \alpha_i - \alpha_j \leq a_j$, in other words
the set
$$
S(a_1, \dots, a_n) :=
\{ (\alpha_1, \dots , \alpha_n ) \, ; \, 0 \leq \alpha_i  \leq \sigma -a_i \}
\bigcap_{1 \leq i < j \leq n} \left (\alpha_j -\alpha_i \geq a_i \quad OR \quad \alpha_i -\alpha_j \geq a_j+1 \right )
$$
is the {\bf singleton} set
$$
\{ (0, a_1, a_1+a_2 , \dots, a_1+ \dots + a_{n-1} ) \} \quad .
$$

{\bf Proof}: 
The condition, for each $1 \leq i<j \leq n$
$$
\alpha_j -\alpha_i \geq a_i \quad OR \quad \alpha_i -\alpha_j \geq a_j+1 
$$
is equivalent to, for each $ 1 \leq i \neq j \leq n$,
$$
\alpha_j -\alpha_i \geq a_i + [i>j] \quad ,
$$
where $[statement]$ is $1$ and $0$ respectively, according to whether the statement is true or false.
Obviously all the $\alpha_i$ are {\bf distinct}, hence there exists a unique permutation $\pi \in S_n$ such that
$$
\alpha_{\pi(1)} < \alpha_{\pi(2)} < \dots < \alpha_{\pi(n)} \quad .
$$
Because of the conditions, we have
$$
\alpha_{\pi(i+1)}-\alpha_{\pi(i)} \geq a_{\pi(i)} + [\pi(i)>\pi(i+1) ]\quad .
$$
Adding up from $i=1$ to $i=n-1$ we have
$$
\alpha_{\pi(n)}-\alpha_{\pi(1)} \geq \sum_{i=1}^{n-1} a_{\pi(i)} + des(\pi) \quad ,
$$
where $des(\pi)$ is the {\it number of descents} of $\pi$  (i.e. the number of $i$, $1 \leq i<n$,  for which $\pi(i)>\pi(i+1)$).
Hence
$$
\alpha_{\pi(n)}-\alpha_{\pi(1)} \geq \sigma-a_{\pi(n)} + des(\pi)  \quad .
$$
{\bf But}, $\alpha_{\pi(n)} \leq  \sigma -a_{\pi(n)}$ and $\alpha_{\pi(1)} \geq 0$ so,
$\alpha_{\pi(n)}-\alpha_{\pi(1)} \leq \sigma -a_{\pi(n)}$, and hence
$$
\sigma-a_{\pi(n)} + des(\pi) \leq \sigma - \alpha_{\pi(n)} \quad .
$$
Hence $des(\pi) \leq 0$. Of course $des(\pi) \geq 0$, hence $des(\pi)=0$, and hence $\pi$ must be
the {\bf identity permutation}: $\pi(i)=i$.

Going back to $\alpha_{\pi(i+1)}-\alpha_{\pi(i)} \geq a_{\pi(i)} + [\pi(i)>\pi(i+1)]$ with $\pi=Identity$, we have
$$
\alpha_1 \geq 0 \quad , \quad 
\alpha_2 -\alpha_1 \geq a_1
\quad , \quad 
\alpha_3 -\alpha_3 \geq a_2
\quad , \quad 
\dots
\quad , \quad 
\alpha_n -\alpha_{n-1} \geq a_{n-1} \quad .
$$
{\bf None} of these inequalities can be strict, or else, adding-them-up would imply that
$\alpha_n > \sigma -a_n$. Hence, the only solution to the above linear-diophantine system of inequalities is
$$
\alpha_1=0 \quad , \quad \alpha_2=a_1 \quad, \quad  \alpha_3=a_1+a_2 \quad , \quad \dots 
\quad , \quad \alpha_n=a_1 + \dots + a_{n-1} \quad . \quad \halmos
$$

{\bf SubFact 5.3}: If $\phi(z)=\prod_{i=0}^{d} (z-q^i)$, and $0 \leq j \leq d$ then
$\phi'(q^j)=\prod_{i=0}^{j-1} (q^j-q^i) \prod_{i=j+1}^{d} (q^j-q^i)$ 
equals $q^{EasyToCompute}(-1)^{AlsoEasyToCompute}$  {\bf times} $(q)_{j} (q)_{d-j}$. \halmos

{\bf SubFact 5.4}: Any evaluation of $F(x_1, \dots, x_n)$ where $x_i=q^{L_i}$ and $L_i$ are affine-linear
expressions in $a_1, \dots, a_n$ can be written in terms of various $(q)_L$ for some affine-linear expressions $L$ 
times  $q^{EasyToCompute}(-1)^{AlsoEasyToCompute}$ \halmos

{\bf Note:} This is better left to a computer, see procedure {\tt EvalFcBG(a,c,r,q,S)} in {\tt qDYSON}.

Plugging-in the unique non-zero point in $S(a_1, \dots, a_n)$, namely $(0, a_1, \dots, a_1 + \dots + a_{n-1})$,
and doing purely-routine
manipulations (better left to the computer), lo-and-behold, we  get what we want, namely
the $q$-multinomial coefficient 
$\frac{(q)_{a_1+a_2+ \dots + a_n}}{(q)_{a_1} (q)_{a_2} \cdots (q)_{a_n}}$.
\halmos. 

(End of proof of Fact 5 (alias the Zeilberger-Bressoud q-Dyson theorem)).

{\bf How to Evaluate Any Other Coefficient}

{\bf Fact 6}:({\bf The Generalized Zeilberger-Bressoud q-Dyson Theorem}) 
Let $\delta=(\delta_1, \dots, \delta_n)$ be
a {\it fixed}, numeric, vector of integers that add-up-to 
$0$. The
coefficient of 
$\prod_{i=1}^{n} x_i^{\delta_i}$ in the $q$-Dyson product
$(qDyson)$ given above is
$$
R_\delta(q,q^{a_1}, \dots, q^{a_n}) \cdot
\frac{(q)_{a_1+a_2+ \dots + a_n}}{(q)_{a_1} (q)_{a_2} \cdots (q)_{a_n}} \quad ,
$$
for some easily-computable (using {\tt qDYSON}) {\it rational function}
$R_{\delta}$. Furthermore, the denominator of $R_{\delta}$ is `nice' (a product of terms of the form $1-q^L$, where
$L$ are affine-linear-combinations in the $a_i$'s).

{\bf Proof}: {\it Now} we apply Fact 4 with
$$
A_i:=\{q^{\alpha_i} \, ; \, 0 \leq 
\alpha_i \leq \sigma -a_i +\delta_i \} \quad .
$$

{\bf SubFact 6.1}: If there exists a pair $i,j$ with  $1 \leq i<j \leq n$ such that
$$
-(a_i-1) \leq \alpha_i - \alpha_j \leq a_j  \quad ,
$$
then $F(q^{\alpha_1}, \dots , q^{\alpha_n})=0$. 

{\bf Proof}: SubFact 6.1 is the same as SubFact 5.1, see the above proof.
\halmos

{\bf SubFact 6.2}: The set of lattice points $(\alpha_1, \dots, \alpha_n)$ with
$0 \leq \alpha_i \leq \sigma - a_i + \delta_i$ such that for {\it every} pair $1 \leq i <j \leq n$
it is not the case that $-(a_i-1) \leq \alpha_i - \alpha_j \leq a_j$, in other words
the set
$$
S_{\delta}(a_1, \dots, a_n) :=
\{ (\alpha_1, \dots , \alpha_n ) \, ; \, 0 \leq \alpha_i  \leq \sigma -a_i+ \delta_i \}
\bigcap_{1 \leq i < j \leq n} \left (\alpha_j -\alpha_i \geq a_i \quad OR \quad \alpha_i -\alpha_j \geq a_j+1 \right )
$$
is a {\bf finite} set, easily constructed by the Maple package {\tt qDYSON}.

{\bf Proof}: 
The condition, for each $1 \leq i<j \leq n$
$$
\alpha_j -\alpha_i \geq a_i \quad OR \quad \alpha_i -\alpha_j \geq a_j+1 
$$
is equivalent to, for each $ 1 \leq i \neq j \leq n$,
$$
\alpha_j -\alpha_i \geq a_i + [i>j] \quad ,
$$
where $[statement]$ is $1$ and $0$ respectively, according to whether the statement is true or false.
Obviously all the $\alpha_i$ are {\bf distinct}, hence there exists a unique permutation $\pi \in S_n$ such that
$$
\alpha_{\pi(1)} < \alpha_{\pi(2)} < \dots < \alpha_{\pi(n)} \quad .
$$
Because of the conditions, we have
$$
\alpha_{\pi(i+1)}-\alpha_{\pi(i)} \geq a_{\pi(i)} + [\pi(i)>\pi(i+1)] \quad .
$$
Adding up from $i=1$ to $i=n-1$, we have
$$
\alpha_{\pi(n)}-\alpha_{\pi(1)} \geq \sum_{i=1}^{n-1} a_{\pi(i)} + des(\pi) \quad ,
$$
where $des(\pi)$ is the {\it number of descents} of $\pi$  (i.e. the number of $i$, $1 \leq i<n$, for which $\pi(i)>\pi(i+1)$).
Hence
$$
\alpha_{\pi(n)}-\alpha_{\pi(1)} \geq \sigma-a_{\pi(n)} + des(\pi)  \quad .
$$
{\bf But}, $\alpha_{\pi(n)} \leq  \sigma -a_{\pi(n)} + \delta_{\pi(n)}$ and $\alpha_{\pi(1)} \geq 0$ so,
$\alpha_{\pi(n)}-\alpha_{\pi(1)} \leq \sigma -a_{\pi(n)} + \delta_{\pi(n)}$, and hence
$$
\sigma-a_{\pi(n)} + des(\pi) \leq \sigma - \alpha_{\pi(n)}+ \delta_{\pi(n)} \quad .
$$
Hence $des(\pi) \leq \delta_{\pi(n)}$. For a given  $\pi \in S_n$ that satisfies this condition, 
we have
$$
\sigma -a_{\pi(n)} + des(\pi) \leq  \sum_{i=1}^{n-1} (\alpha_{\pi(i+1)} - \alpha_{\pi(i)} ) \leq \sigma -a_{\pi(n)} + \delta_{\pi(n)} \quad .
$$
Since
$$
\alpha_{\pi(i+1)}-\alpha_{\pi(i)} \geq a_{\pi(i)} + [\pi(i)>\pi(i+1)] \quad .
$$
we can write, for $i=1, \dots, n-1$, and for some integers $m_i \geq 0$ ($1 \leq i \leq n$):
$\alpha_{\pi(1)}=\alpha_1$, and for $1 \leq i \leq n-1$:
$$
\alpha_{\pi(i+1)}-\alpha_{\pi(i)} \, = \, a_{\pi(i)} + [\pi(i)>\pi(i+1)]+ m_{i+1} \quad .
$$
Summing from $i=1$ to $i=n-1$ gives
$$
0 \leq \sum_{i=1}^{n} m_i \leq \delta_{\pi(n)} -des(\pi) \quad .
$$
Of course, there are only finitely-many such $\{(m_1, \dots, m_n)\}$. 
So, for {\it each} permutation $\pi$ obeying $des(\pi) \leq \delta_{\pi(n)}$ and for {\it each}
vector $(m_1, m_2, \dots, m_n)$ of non-negative integers whose sum is $\leq \delta_{\pi(n)}- des(\pi)$,
we have  a member of $S_{\delta}(a_1, \dots, a_n)$,
$$
\alpha_{\pi(i)}=\sum_{r=1}^{i} m_r +  \sum_{r=1}^{i-1} ( a_{\pi(r)}+[\pi(r)>\pi(r+1)] ) \quad .
$$
Note that $a_i$ are {\it symbolic}, and the sets of feasible $\pi$ and $(m_1, \dots, m_n)$ {\it only}
depend to $\delta$ {\bf not} on $(a_1, \dots, a_n)$.
It is immediately seen that these points satisfy the condition for membership in $S_\delta$. \halmos

{\bf SubFact 6.3}: For a given feasible permutation $\pi$ and 
feasible vector $m=(m_1 \dots , m_n )$,
and with $A_i$ as above, the summand of Fact 4 is a {\bf simple} (factored)
rational function of $(q, q^{a_1}, \dots, q^{a_n})$ {\bf times} the
$q$-multinomial coefficient
$\frac{(q)_{a_1+a_2+ \dots + a_n}}{(q)_{a_1} (q)_{a_2} \cdots (q)_{a_n}}$. 

{\bf Proof}: Routine (and programmed into {\tt qDYSON}). \halmos

Adding up these {\it finitely} many contributions concludes the proof
of Fact 6. \halmos

{\bf Remark}: One can get much smaller sets of evaluation-points 
$S_\delta$, by
shifting the $A_i$'s by (positive or negative) $c_i$, in other words
consider
$$
A_i:=\{q^{\alpha_i} \, ; \, c_i \leq \alpha_i \leq \sigma -a_i +\delta_i+c_i  \} \quad .
$$
Of course, one should get the {\it same} output, regardless of the $c$'s, but
for the sake of efficiency it would be nice to make $S_\delta$ as small as possible.
The Maple package {\bf qDYSON} has a procedure {\tt BestShift(d)}, that
finds the optimal shift.

{\bf Another Remark}: Fact 6 is only valid for
{\it numeric} (specific) $\delta$, and each numeric, specific, number of variables $n$. 
There is no closed form formula for the general coefficient of 
$\prod_{i=1}^{n} x_i^{\delta_i}$ of the $q$-Dyson product
where the $\delta$, as well as $a=(a_1, \dots, a_n)$,
are symbolic. For any specific $n$, it follows from
WZ theory[WZ], that this quantity is {\it holonomic},
i.e. there exist `pure' linear recurrences (in {\bf each} of $a_i$) with coefficients
that are polynomials in $(q^{a_1}, \dots, q^{a_n}, q^{\delta_1}, \dots, q^{\delta_n})$,
but these are already fairly complicated for $n=3$, and the
orders get larger and larger with larger $n$.

The miracle of $q$-Dyson is that for the  constant term,
these recurrences are always {\it first-order}
(that is what it means to be {\it closed-form}), for {\bf any} number
of variables, $n$, and the generalized Zeilberger-Bressoud
(Fact 6) extends this miracle to {\it any},
specific, other coefficient.

{\bf Yet another remark}: Even though Fact 6 is only valid, in general, for
specific $n$ and specific $\delta$, it sometimes happens that
if you take a specific, numeric, $\delta=(\delta_1, \dots, \delta_n)$
then the coefficient of $\delta_r:=(\delta_1, \dots, \delta_n, 0( r \,\,\, times))$
in the $q$-Dyson product in $n+r$ variables can be expressed `uniformly', since
the size of $S_{\delta_r}$ (shifted by a judicious shift)
remains the same.
Of course, this is the case with the original $q$-Dyson, with $\delta=(0, \dots, 0)$, 
where $S_{\delta}$ is always a singleton, leading to a general statement valid for {\it all} $n$
(i.e. symbolic $n$).
This is also the case with the conjectures of
Drew Sills[S2] proved 
by Lun Lv, Guoce Xin and Yue Zhou [LXZ], by extending the method of [GX].
In [LXZ] there is also a beautiful, much more general theorem. 
We are sure that their result can be reproved, quicker, using the
K\'aroly-Nagy approach, as extended in our present article, but we leave this
to the interested reader.

{\bf The Maple package qDYSON}

Everything (and more) is implemented in the Maple package
{\tt qDYSON}. 
Its only limitation is that the number
of variables, $n$, is {\it numeric}, not symbolic,
but, as noted above, often, by running it for
$n \leq 5$, one can deduce, and easily
translate for general $n$, the proofs given by the
package for specific $n$.

The main procedure is: \quad `{\tt Gyula(z,d,q);}'.
It inputs a variable-name, $z$,  a numeric list of intgers (of length $n$, say),
$d=[d_1, \dots, d_n]$ whose sum is $0$, 
as well as a variable-name, $q$.
It outputs
the rational function $R_d(q, q^{a_1}, \dots, q^{a_n})$
promised by Fact $6$, where, for the sake of clarity, $q^{a_i}$ is replaced by $z_i$. 
For example,

`{\tt Gyula(z,[2,-2,0,0],q);}' yields 
$$
( -{z_{{4}}}^{2}{z_{{2}}}^{2}{q}^{3}z_{{3}}z_{{1}}+{z_{{4}}}^{3}{z_{{2}}}^{2}{q}^{3}z_{{3}}z_{{1}}+z_{{3}}{z_{{2}}}^{2}z_{{4}}{q}^{3}-{z_{{2}}}^{2
}{q}^{3}{z_{{3}}}^{2}{z_{{4}}}^{2}z_{{1}}+{z_{{2}}}^{2}{q}^{3}{z_{{3}}}^{3}{z_{{4}}}^{2}z_{{1}}-{z_{{2}}}^{2}z_{{4}}{q}^{2}z_{{1}}z_{{3}}-
$$
$$
z_{{2}}{q}^{2}z
_{{3}}-{z_{{4}}}^{2}z_{{2}}{q}^{2}+z_{{2}}z_{{4}}{q}^{2}z_{{1}}z_{{3}}-z_{{2}}{z_{{3}}}^{2}z_{{4}}{q}^{2}+z_{{3}}z_{{2}}z_{{4}}{q}^{2}-{z_{{3}}}^{2}{z_{{
4}}}^{2}z_{{2}}{q}^{2}z_{{1}}+z_{{2}}z_{{4}}qz_{{1}}+z_{{1}}z_{{3}}z_{{2}}q+z_{{3}}z_{{4}}q-z_{{1}} )   \cdot
$$
$$
\frac{\left( 1-z_{{1}} \right)}{
 \left( z_{{4}}z_{{2}}q-1 \right)  \left( z_{{3}}z_{{2}}q-1 \right)  \left( z_{{2}}z_{{3}}z_{{4}}q-1 \right)  \left( z_{{3}}z_{{2}}z_{{4}}{q}^{2}-1
 \right)} \quad .
$$
As $d$ and/or $n$ gets larger, the size of $S_{d}$ gets
larger, and Maple takes a long time to
bring everything under a common denominator 
(using the command {\tt normal}), so for these, a much
faster alternative is the unsimplified version, `{\tt Zoltan (z,d,q);}' ,
that in some sense is better, since it displays the output
as a sum of simple rational functions. For example, `{\tt Zoltan(z,[2,-2,0,0],q);}` yields
$$
-{\frac { \left( -1+z_{{1}} \right)  \left( -z_{{1}}+q \right) }{ \left( qz_{{2}}z_{{3}}z_{{4}}-1 \right)  \left( z_{{3}}z_{{2}}z_{{4}}{q}^{2}-1 \right) 
}}+{\frac {z_{{2}}{q}^{2} \left( -1+z_{{1}} \right)  \left( z_{{4}}z_{{1}}-1 \right)  \left( -1+z_{{3}} \right)  \left( z_{{2}}qz_{{3}}z_{{4}}z_{{1}}-1
 \right) }{ \left( -1+z_{{2}}q \right)  \left( z_{{4}}z_{{2}}q-1 \right)  \left( z_{{2}}z_{{3}}q-1 \right)  \left( z_{{3}}z_{{2}}z_{{4}}{q}^{2}-1
 \right) }}
$$
$$
+{\frac {z_{{4}}z_{{2}}{q}^{2} \left( -1+z_{{4}} \right)  \left( -1+z_{{1}} \right)  \left( z_{{3}}z_{{1}}-1 \right)  \left( z_{{2}}qz_{{3}}z_
{{4}}z_{{1}}-1 \right) }{ \left( z_{{4}}z_{{2}}q-1 \right)  \left( z_{{2}}z_{{3}}q-1 \right)  \left( qz_{{2}}z_{{3}}z_{{4}}-1 \right)  \left( z_{{3}}z_{{
2}}z_{{4}}{q}^{2}-1 \right) }}-{\frac {q \left( -1+z_{{3}}z_{{4}} \right)  \left( -1+z_{{1}} \right)  \left( z_{{2}}z_{{4}}qz_{{1}}-1 \right)  \left( z_{
{2}}qz_{{3}}z_{{4}}z_{{1}}-1 \right) }{ \left( -1+z_{{2}}q \right)  \left( z_{{4}}z_{{2}}q-1 \right)  \left( qz_{{2}}z_{{3}}z_{{4}}-1 \right)  \left( z_{
{3}}z_{{2}}z_{{4}}{q}^{2}-1 \right) }} \quad .
$$

Other important procedures are 
{\tt GyulaTh} , {\tt ZoltanTh}, that are verbose
forms of the above, and  {\tt Sefer},  that
outputs full articles. Some examples are given in the front of this article

{\eighttt http://www.math.rutgers.edu/\~{}zeilberg/mamarim/mamarimhtml/qdyson.html} .

Also of interest is procedure `{\tt BestShift(d);}',  already mentioned above, that finds the
best shift to make $S_d$ as small as possible.
There are also plenty of checking procedures, that make
sure that everything is correct! Please consult
the on-line help 
(gotten by typing {\tt ezra();}, {\tt ezra1();}, {\tt ezraS(); }, and {\tt ezraC();}).

{\bf A very brief history}

The original Dyson conjecture appeared in 1962, 
in a very important paper[Dy] (that according to {\tt google scholar}  (viewed Aug. 13, 2013)
was cited 1483 times), and (whose sequel) lead to, {\it inter alia},
intriguing connections to the Riemann Zeta function
discovered by Hugh Montgomery, and extended by Andrew
Odlyzko, Jonathan Keating, Nina Snaith and others. 
The Dyson conjecture was proved
shortly after by Jack Gunson[Gu] and 
Kenneth G. Wilson(1936-2013, Physics Nobel 1982)[W], who later on went on  to revolutionize
physics by creating {\it renormalization group theory}.
The {\it proof from the book}
of the original Dyson conjecture 
(using  Fact 3 with $d=n-1$) was given[Go] in 1970 by 
Bayesian pioneer (and collaborator of Alan Turing at Bletchley Park) 
Jack Good(1916-2009). 
In 1982, one of us (Zeilberger) found a longer, but equally nice, {\it combinatorial proof}[Z].
A couple of years later, Dave Bressoud and Doron Zeilberger succeeded in $q$-ifing this proof,
thereby giving the {\bf first} proof of the $q$-analog, conjectured, in 1975, by George Andrews[An].
A shorter proof of the Zeilberger-Bressoud $q$-Dyson
theorem was given by Ira Gessel and Guoce Xin [GX],
and as already noted, the {\it proof from the book}
was given by K\'arolyi and Nagy[KN], 
that formed the inspiration
for the present article.
Other far-reaching applications of their method
are given  by Gyula K\'arolyi in collaboration with
Alain Lascoux and Ole Warnaar [KLW].

The problem of computing other coefficients, besides
the constant term, for the original Dyson product,
was launched by Sills and Zeilberger[SZ], followed by
Sills' more general article [S1], that was $q$-analogized in [S2], where he
conjectured interesting `uniform' expressions for a few other
coefficients of the $q$-Dyson, proved,
and vastly generalized, in the article [LXZ], already 
mentioned above.

{\bf Yet another approach}

While we {\it love} the K\'arolyi-Nagy proof, as
extended in this article, let us end by remarking
that the {\it original} approach of [ZB] could also be
used to prove Fact 6, and even implemented on
a computer. Now the multi-tournaments 
of [Z] and [ZB] have different
score vectors, and unlike the original case,
where the `good guys' can be completely translated into
words, with no `left-overs', now, we do get `minimal left-overs'.
It is easy to see
that, for specific $\delta$, and specific $n$,
there are only finitely many of them, and each
of them could be $q$-counted. But, since the
K\'arolyi-Nagy approach is {\it so} efficient, there is
little motivation to extend the original
Zeilberger-Bressoud approach.

{\bf References}

[Al] N. Alon, {\it Combinatorial Nullstellensatz}, Combin. Probab. Comput. {\bf 8}(1999) 7-29. \hfill\break
{\tt http://www.tau.ac.il/\~{}nogaa/PDFS/null2.pdf}

[An] G. E. Andrews, {\it Problems and prospects for basic hypergeometric functions}, 
in: ``Theory and Application of Special Functions'', (R. A. Askey, ed.), Academic Press, New York (1975), 191-224.

[Dy] F. Dyson, {\it Statistical theory of energy levels of complex systems}, J. Math. Phys. {\bf 3}(1962), 140-156.

[GX] I. M. Gessel and G. Xin, {\it A short proof of the Zeilberger–Bressoud $q$-Dyson theorem}, Proc. Amer. Math. Soc. {\bf 134}(2006), 2179-2187. \hfill\break
{\tt http://arxiv.org/abs/math/0412339}

[Go] I. J. Good, {\it Short proof of a conjecture by Dyson}, J. Math. Phys. {\bf 11}(1970) 1884.

[Gu] J. Gunson, {\it Proof of a conjecture of Dyson in the statistical theory of energy levels}, J. Math. Phys. {\bf 3}(1962) 752-753.

[KLW] Gy. K\'arolyi, A. Lascoux, and S. O. Warnaar, {\it Constant term identities and Poincar\'e polynomials}, \hfill\break
{\tt http://arxiv.org/abs/1209.0855} (also to appear in Trans. Amer. Math. Soc.).

[KN] Gy. K\'arolyi and Z.L. Nagy, {\it A simple proof of the Zeilberger-Bressoud $q$-Dyson theorem}, \hfill\break
{\tt http://arxiv.org/abs/1211.6484} (also to appear in Proc. Amer. Math. Soc.) .

[KP] R.N. Karasev and F.V. Petrov, {\it Partitions of non-zero elements of a finite field into pairs},
Israel J. of Math., {\bf 192} (2012), 143-156. arxiv.org/abs/1005.1177

[L] M. Laso\'n, {\it A generalization of Combinatorial Nullstellensatz}, Electron. J. Combin. {\bf 17}(2010) \#N32, 6 pages. 
{\tt http://arxiv.org/abs/1302.4647}

[LXZ] L. Lv, G. Xin and Y. Zhou, {\it A Family of $q$-Dyson Style Constant Term Identities}, J. Combin. Theory Ser. A {\bf 116} (2009), 12-29.
{\tt  http://arxiv.org/abs/0706.1009}

[S1] A. V. Sills, {\it Disturbing the Dyson conjecture, in a generally GOOD way}, J. Combin. Theory Ser. A {\bf 113} (2006), 1368-1380.
\hfill\break
{\tt http://math.georgiasouthern.edu/\~{}asills/Dyson2/dyson2.pdf}

[S2] A. V. Sills, {\it Disturbing the $q$-Dyson conjecture}, in ``Tapas in Experimental Mathematics'', (T. Amdeberhan and
V. Moll, eds.), Contemporary Mathematics {\bf 457}, 265-272.
\hfill\break
{\tt http://math.georgiasouthern.edu/\~{}asills/qDysonDist/DysonMethodFinal.pdf}

[SZ] A.V. Sills and D. Zeilberger,
{\it Disturbing the Dyson Conjecture (in a GOOD Way)},
Experimental Mathematics {\bf 15}(2006), 187-191.
\hfill\break
{\tt http://www.math.rutgers.edu/\~{}zeilberg/mamarim/mamarimPDF/gooddyson.pdf}

[W] K.G. Wilson, {\it Proof of a conjecture by Dyson, J. Math}. Phys. {\bf 3}(1962), 1040-1043.

[WZ] H. Wilf and D. Zeilberger, 
{\it An algorithmic proof theory for hypergeometric (ordinary and ``q'') multisum/integral identities},
Invent. Math. {\bf 108}(1992), 575-633.
\hfill\break
{\tt http://www.math.rutgers.edu/\~{}zeilberg/mamarim/mamarimPDF/multiwz.pdf}

[Z]  D. Zeilberger, {\it A combinatorial proof of Dyson's conjecture}, Discrete Math. {\bf 41}(1982), 317-321.
\hfill\break
{\tt http://www.math.rutgers.edu/\~{}zeilberg/mamarimY/dysonconj.pdf}

[ZB] D. Zeilberger and D. Bressoud, {\it A proof of Andrews's $q$-Dyson conjecture}, Discrete Math. {\bf 54}(1985), 201-224.
\hfill\break
{\tt http://www.math.rutgers.edu/\~{}zeilberg/mamarimY/ZeilbergerBressoudTheorem.pdf}

\end